\documentclass[12pt]{amsart}

\usepackage{amsmath,xspace,amssymb,mathrsfs}
\usepackage{color}

\usepackage{hyperref}
\usepackage{enumerate}
\usepackage{a4wide}
\usepackage{color}

\input xy
\xyoption{all}
\xyoption{2cell}
\UseAllTwocells
\CompileMatrices

\newcommand{\Spec}{\operatorname{Spec}}
\renewcommand{\phi}{\varphi}

\newcommand{\Ker}{\operatorname{Ker}}

\newcommand{\Max}{\operatorname{Max}}

\newcommand{\Min}{\operatorname{Min}}

\newcommand{\Clop}{\operatorname{Clop}}
\newcommand{\Fin}{\operatorname{Fin}}

\newcommand{\Sp}{\operatorname{Sp}}

\newtheorem{proposition}{Proposition}[section]
\newtheorem{lemma}[proposition]{Lemma}

\newtheorem{corollary}[proposition]{Corollary}
\newtheorem{theorem}[proposition]{Theorem}

\theoremstyle{definition}

\newtheorem{example}[proposition]{Example}

\newtheorem{remark}[proposition]{Remark}

\begin{document}

\title[Lifting, orthogonality and finiteness of idempotents]{Structural results on lifting, orthogonality and finiteness of idempotents}

\author[A. Tarizadeh, P.K. Sharma]{Abolfazl Tarizadeh, Pramod K. Sharma}
\address{Department of Mathematics, Faculty of Basic Sciences, University of Maragheh \\
P. O. Box 55136-553, Maragheh, Iran.}
\email{ebulfez1978@gmail.com}
\address{Retired Professor, School of Mathematics, D.A.V.V. \\ Indore-452001, India.}
\email{pksharma1944@yahoo.com}

\date{}

\begin{abstract}
In this paper, using the canonical correspondence between the idempotents and clopens, we obtain several new results on lifting idempotents. The Zariski clopens of the maximal spectrum are precisely determined, then as an application, lifting idempotents modulo the Jacobson radical is characterized. Lifting idempotents modulo an arbitrary ideal is also characterized in terms of certain connected sets related to that ideal. Then as an application, we obtain that the sum of a lifting ideal and a regular ideal is a lifting ideal. We prove that lifting idempotents preserves the orthogonality in countable cases. The lifting property of an arbitrary morphism of rings is characterized. As another major result, it is proved that the number of idempotents of a ring $R$ is finite if and only if it is of the form $2^{\kappa}$ where $\kappa$ is the cardinal of the connected components of $\Spec(R)$. Finally, it is proved that the primitive idempotents of a zero dimensional ring are in 1-1 correspondence with the isolated points of its prime spectrum. These results either generalize or improve several important results in the literature.
\end{abstract}

\maketitle

\section{Introduction}

Many structural aspects of rings, decompositions of modules (especially Peirce and Wedderburn decompositions), representation theory of associative algebras, projective covers and some homological properties of rings are tied up with idempotents and their lifting properties. For instance, it is well known that a ring $R$ is a clean ring (i.e., each $f\in R$ can be written as $f=e+g$ where $e\in R$ is an idempotent and $g\in R$ is invertible) if and only if its idempotents can be lifted modulo each ideal (see \cite{Nicholson}). In a recent work \cite[Theorem 8.5]{Tarizadeh}, we proved the dual of this fundamental result which states the idempotents of a reduced ring $R$ can be lifted along each localization if and only if it is a purified ring (i.e., every pair of distinct minimal prime ideals of $R$ can be separated by an idempotent). As another application, lifting idempotents are useful in understanding the structure of p.p. rings and their generalizations (for further information see e.g. \cite{A. Tarizadeh Racsam}). Lifting idempotents of rings (especially non-commutative rings) modulo an ideal have been periodically studied in the literature over the years, see e.g. \cite{Tarizadeh}-\cite{Zhu}. \\

This paper mainly focuses on lifting, orthogonality and investigating the cardinality of idempotents. Several results of this paper are based on using the folklore correspondence between the idempotents and clopens. Indeed, the following fundamental result has numerous applications in commutative algebra and algebraic geometry. It also plays a major role and unifying framework in this paper.

\begin{theorem}\label{Theorem II} The map $e\rightsquigarrow D(e)$ is a one-to-one correspondence from the set of idempotents of a ring $R$ onto the set of clopen $($both open and closed$)$ subsets of $\Spec(R)$.
\end{theorem}

{\bf Proof.} See \cite[00EE]{Johan}. $\Box$ \\

By $\mathcal{B}(R)=\{e\in R: e=e^{2}\}$ we mean the set of all idempotents of a given ring $R$. Remember that this set with the operations $e\oplus e'=e+e'-2ee'$ as the addition and $e.e'=ee'$ as the multiplication is a ring. We call $\mathcal{B}(R)$ the Boolean ring of $R$. The set of clopens of a topological space $X$ is denoted by $\Clop(X)$. It has a ring structure. In fact, $\Clop(X)$ is a subring of the power set ring $\mathcal{P}(X)$. If $\phi:X\rightarrow Y$ is a continuous map of topological spaces, then the induced map $\Clop(\phi):\Clop(Y)\rightarrow\Clop(X)$ given by $U\rightsquigarrow\phi^{-1}(U)$ is a morphism of rings. In fact,
both formations $\mathcal{B}(-)$ and $\Clop(-)$ are covariant and contravariant functors, respectively. For more information see \cite[\S3]{Tarizadeh 2}. It is worth mentioning that the correspondence of Theorem \ref{Theorem II} is in fact an isomorphism of rings $\mathcal{B}(R)\simeq\Clop\big(\Spec(R)\big)$. This observation, in particular, generalizes Stone representation theorem to the arbitrary rings (see \cite[Theorem 3.1]{Tarizadeh 2}). \\

The following is a brief outline of the paper. In Theorem \ref{Corollary iv nice}, the Zariski clopens of the maximal spectrum are precisely determined. Then as an application, in Corollary \ref{Theorem v new}, we obtain the canonical isomorphism of rings $\mathcal{B}(R/\mathfrak{J})\simeq\Clop\big(\Max(R)\big)$.
Moreover lifting idempotents modulo the Jacobson radical is characterized (see Corollary \ref{Corollary v Jacobson}). Theorem \ref{Theorem vi new} characterizes lifting idempotents modulo an arbitrary ideal, especially in the part (v), the lifting property of an ideal is characterized in terms of certain connected sets related to that ideal. This is a very useful criterion to check the lifting property of an ideal. In particular, we obtain that the sum of a lifting ideal and a regular ideal is a lifting ideal (see Corollary \ref{Corollary vi regular}). Corollary \ref{Theorem I} and Proposition \ref{Proposition II} provide further criteria to check the lifting property of an ideal. In Corollary \ref{Corollary viii along locz}, lifting idempotents along the localizations is investigated (also see \cite[Theorem 8.5(vi)]{Tarizadeh}). We also provide new proofs for some major results in the literature on the lifting idempotents (see Propositions \ref{Theorem III} and \ref{Theorem viii complete}). \\

Section 4 investigates the cardinality of idempotents. Especially we prove that for any ring $R$, if $|R|\leqslant|\Spec(R)|$ then $\Spec(R)$ is infinite (see Proposition \ref{proposition vi zero dim}). If moreover $R$ is a Boolean ring, then the converse holds (see Theorem \ref{Theorem Deligne-Tarizadeh}). Then as an application, in Theorem \ref{Corollary ix main}, it is proved that the number of idempotents of a ring $R$ is finite if and only if it equals $2^{\kappa}$ where $\kappa$ is the cardinal of the connected components of $\Spec(R)$. \\

Section 5 investigates the orthogonal and primitive idempotents. Especially in Theorem \ref{Theorem  vii orthogonality}, we prove that a countable set of orthogonal idempotents can be lifted orthogonally by every idempotent lifting morphism of rings. The lifting property of an arbitrary morphism of rings is characterized (see Lemma \ref{Lemma referee com} and Theorem \ref{Theorem lift h0}). It is also proved that the primitive idempotents of a zero dimensional ring are in 1-1 correspondence with the isolated points of its prime spectrum (see Theorem \ref{Theorem x isolated points}). Its proof uses the following result which is interesting in its own right: the structure of prime ideals of a zero dimensional ring $R$ and its Boolean ring $\mathcal{B}(R)$ are the same. Theorem \ref{Theorem x isolated points}, in particular, yields that there are infinite Boolean rings whose primitive idempotents are finite. Hence, contrary to the connected components, the cardinal of primitive idempotents is not an invariant to compute the cardinal of all idempotents. Theorem \ref{Theorem xii TD} is the further main result of this section. \\

Finally, Appendix (\S6) concludes with a new and short proof of the classical Chinese remainder theorem. Needless to say that this theorem has numerous applications in mathematics, (for a most recent application of it see  \cite{Kowalski}). \\

Throughout the paper we give various non trivial examples which illuminate the involved abstract structures.

\section{Preliminaries}

We collect in this section some basic background for the reader's convenience. \\

In this paper, all rings are commutative. If $f$ is a member of a ring $R$, then $D(f)=\{\mathfrak{p}\in\Spec(R): f\notin\mathfrak{p}\}$ and $V(f)=\Spec(R)\setminus D(f)$. The nilradical and Jacobson radical of a ring $R$ is denoted by $\mathfrak{N}$ and $\mathfrak{J}$, respectively. If $\mathfrak{p}$ is prime ideal of a ring $R$, then the canonical ring map $R\rightarrow R_{\mathfrak{p}}$ is denoted by $\pi_{\mathfrak{p}}$. If $\phi:R\rightarrow R'$ is a morphism of rings then the induced map $\Spec(\phi):\Spec(R')\rightarrow\Spec(R)$ given by $\mathfrak{p}\rightsquigarrow\phi^{-1}(\mathfrak{p})$ is also denoted by $\phi^{\ast}$. \\

We say that a morphism of rings $\phi:R\rightarrow R'$ lifts idempotents if whenever $e'\in R'$ is an idempotent, then there exists an idempotent $e\in R$ such that $\phi(e)=e'$. If moreover such an idempotent $e$ is unique then we say that $\phi$ strongly lifts idempotents. Clearly a morphism of rings $\phi: R\rightarrow R'$ strongly lifts idempotents if and only if the ring map $\mathcal{B}(\phi):\mathcal{B}(R)\rightarrow\mathcal{B}(R')$ given by $e\rightsquigarrow\phi(e)$ is an isomorphism. Lifting idempotents is stable under the composition of ring maps. Let $\psi:R'\rightarrow R''$ be a second morphism of rings. If $\psi\circ\phi$ lifts idempotents, then $\psi$ is as well. If moreover $\psi$ is injective, then $\phi$ also lifts idempotents. For any ring $R$, the canonical injections $R\rightarrow R[x]\rightarrow R[[x]]$ lift idempotents. In other words, if $f=\sum\limits_{i\geqslant0}r_{i}x^{i}\in R[[x]]$ is an idempotent, then $f=r_{0}$. For its proof see \cite[Lemma 2.3]{A. Tarizadeh Racsam}.
In the same vein, if $R$ is a $\mathbb{Z}$-graded ring then the ring extension $R_{0}\subseteq R$ lifts idempotents. In other words, every nonzero idempotent of a $\mathbb{Z}$-graded ring is homogeneous of degree zero. For its proof see \cite[Theorem 1]{Kirby} or \cite[Theorem 6.4]{Tarizadeh-Oinert}. \\

Let $I$ be an ideal of a ring $R$. If the canonical ring map $R\rightarrow R/I$ lifts idempotents, then it is also called that $I$ is a lifting ideal of $R$ (or, the idempotents of $R$ can be lifted modulo $I$). Thus $I$ is a lifting ideal of $R$ if and only if $f-f^{2}\in I$ for some $f\in R$, then there exists an idempotent $e\in R$ such that $f-e\in I$. For simplicity, we often write ``lifting ideal'' instead of ``lifting ideal of $R$''. If $R/I$ has no nontrivial idempotents, then $I$ is a lifting ideal. But the converse does not hold. For example, the zero ideal is always a lifting ideal, but there are many rings with nontrivial idempotents, (for instance $\mathbb{Z}/6\mathbb{Z}$ has nontrivial idempotents 3 and 4). Let $S$ be a multiplicative subset of a ring $R$. If the canonical ring map $R\rightarrow S^{-1}R$ lifts idempotents then we say that the idempotents of $R$ can be lifted along the localization $S^{-1}R$. \\

One can see that $m\mathbb{Z}$ is a lifting ideal of $\mathbb{Z}$ if and only if $m=0,1, p^{d}$ where $p$ is a prime number and $d\geqslant1$. Indeed, if $m\mathbb{Z}$ is a lifting ideal with $m\geqslant2$, then we show that $m$ is of the form $p^{d}$. If not, then by the Chinese remainder theorem, the ring $\mathbb{Z}/m\mathbb{Z}$ contains a nontrivial idempotent $a+m\mathbb{Z}$. Since $m\mathbb{Z}$ is a lifting ideal, thus there exists an idempotent $e\in\{0,1\}$ in $\mathbb{Z}$ such that $a-e\in m\mathbb{Z}$. This yields that either $a\in m\mathbb{Z}$ or $a-1\in m\mathbb{Z}$. But this is a contradiction, because $a+m\mathbb{Z}$ is a nontrivial idempotent. Conversely, if $m\in\{0,1\}$ then the canonical ring map $\mathbb{Z}\rightarrow\mathbb{Z}/m\mathbb{Z}$ is either the identity map or the zero map which both maps obviously lift idempotents. If $m=p^{d}$ then $\mathbb{Z}/p^{d}\mathbb{Z}$ is a local ring and so its idempotents are trivial. This completes the argument. Lifting property of ideals need not pass to neither intersections nor products. For example, if $p$ and $q$ are distinct prime numbers then $p\mathbb{Z}$ and $q\mathbb{Z}$ are lifting ideals of $\mathbb{Z}$, but $p\mathbb{Z}\cap q\mathbb{Z}=pq\mathbb{Z}$ is not a lifting ideal. \\

If $X$ is a set, then by $\Fin(X)$ we mean the set of all finite subsets of $X$. The number of elements (cardinal) of a set $X$ is denoted by $|X|$. By $\pi_{0}(X)$ we mean the space of connected components of a topological space $X$ which carries the quotient topology. By a regular ideal of a ring $R$ we mean an ideal of $R$ which is generated by a set of idempotents of $R$. Note that the term ``regular ideal'' may have different meanings in the literature. Each maximal element of the set of proper and regular ideals of a ring $R$ is called a max-regular ideal of $R$. The set of max-regular ideals of $R$ is called the \emph{Pierce spectrum} of $R$ and is denoted by $\Sp(R)$. It is a compact and totally disconnected space whose base opens are precisely of the form $U_{e}=\{M\in\Sp(R): e\notin M\}$ where $e\in R$ is an idempotent. For more details on this topic see \cite{Tarizadeh}, \cite{Tarizadeh 5} and \cite{Tarizadeh 2}. Finally, we have the canonical isomorphisms of topological spaces $\Spec\big(\mathcal{B}(R)\big)\simeq\Sp(R)
\simeq\pi_{0}\big(\Spec(R)\big)$.
For the first isomorphism see \cite[Theorem 4.1]{Tarizadeh 2}. The second isomorphism is given by $M\rightsquigarrow V(M)$. \\

In dealing with Theorem \ref{Theorem II}, it is important to notice that if $D(f)$ is a clopen of $\Spec(R)$ for some $f\in R$, then there exists an idempotent $e\in R$ such that $D(f)=D(e)$. But we cannot deduce that $f=e$. For example, $D(-1)=\Spec(\mathbb{Z})=D(1)$ but $1\neq-1$.

\section{Lifting idempotents modulo an ideal}

As the first application of Theorem \ref{Theorem II} we obtain the following characterization.

\begin{lemma}\label{Lemma referee com} A morphism of rings $\phi: R\rightarrow R'$ lifts idempotents if and only if the induced map $\Clop(\phi^{\ast}):\Clop\big(\Spec(R)\big)
\rightarrow\Clop\big(\Spec(R')\big)$ is surjective.
\end{lemma}

{\bf Proof.} It easily follows from Theorem \ref{Theorem II} and the fact that if $e\in R$ is an idempotent then $\Clop(\phi^{\ast})\big(D(e)\big)=
D\big(\phi(e)\big)$. $\Box$

\begin{corollary}\label{Theorem I} Let $I$ be an ideal of a ring $R$. Then $I$ is a lifting ideal of $R$ if and only if $\sqrt{I}$ is as well.
\end{corollary}

{\bf Proof.} The canonical ring map $\eta:R\rightarrow R/\sqrt{I}$ factors as $\eta=\mu\circ\pi$ where $\pi:R\rightarrow R/I$ and $\mu:R/I\rightarrow R/\sqrt{I}$ are the canonical ring maps. It follows that $\Clop(\eta^{\ast})=\Clop(\mu^{\ast})\circ\Clop(\pi^{\ast})$. Note that $\mu^{\ast}$ is a homeomorphism and so $\Clop(\mu^{\ast})$ is bijective. Now using these, then the assertion easily follows from Lemma \ref{Lemma referee com}. $\Box$

\begin{corollary} Let $I$ and $J$ be two ideals of a ring $R$ such that $\sqrt{I}=\sqrt{J}$. Then $I$ is a lifting ideal if and only if $J$ is so.
\end{corollary}
{\bf Proof.} It is an immediate consequence of Corollary \ref{Theorem I}. $\Box$

\begin{corollary}\label{Corollary I} Every ideal contained in the nilradical is a lifting ideal.
\end{corollary}

{\bf Proof.} Let $I$ be an ideal of a ring $R$ such that $I\subseteq\mathfrak{N}$. So $\sqrt{I}=\sqrt{0}$. Then apply Corollary \ref{Theorem I}. $\Box$

\begin{lemma}\label{prop s.lift i2} A morphism of rings $\phi:R\rightarrow R'$ strongly lifts idempotents if and only if the map $\Clop(\phi^{\ast})$ is bijective.
\end{lemma}

{\bf Proof.} The following diagram of Boolean rings is commutative: $$\xymatrix{\mathcal{B}(R)
\ar[r]^{\mathcal{B}(\phi)}\ar[d]^{\simeq}&\mathcal{B}(R')
\ar[d]^{\simeq}\\
\Clop\big(\Spec(R)\big)\ar[r]^{\Clop(\phi^{\ast})}&\Clop\big(\Spec (R')\big)}$$
where the vertical arrows are the canonical isomorphisms. Then $\mathcal{B}(\phi)$ is bijective if and only if $\Clop(\phi^{\ast})$ is as well. $\Box$

\begin{corollary} Let $\phi:R\rightarrow R'$ be a morphism of rings. If $\phi^{\ast}:\Spec(R')\rightarrow\Spec(R)$ is a homeomorphism, then $\phi$ strongly lifts idempotents.
\end{corollary}

{\bf Proof.} By hypothesis, $\Clop(\phi^{\ast})$ is bijective. Then apply Lemma \ref{prop s.lift i2}. $\Box$

\begin{corollary} Let $I\subseteq J\subseteq\sqrt{I}$ be ideals of a ring $R$. Then the canonical ring map $\pi:R/I\rightarrow R/J$ strongly lifts idempotents.
\end{corollary}

{\bf Proof.} Clearly $V(I)=V(J)$. Thus $\pi^{\ast}$ is a homeomorphism and so $\Clop(\pi^{\ast})$ is bijective. Then apply Lemma \ref{prop s.lift i2}. $\Box$

\begin{lemma}\label{Lemma II} Let $\mathfrak{p}$ be a prime ideal of a ring $R$ and $e,e'\in R$ be idempotents. If $e-e'\in\mathfrak{p}$ then $e-e'\in\Ker\pi_{\mathfrak{p}}$. In particular, if $e-e'\in\mathfrak{J}$ then $e=e'$.
\end{lemma}

{\bf Proof.} Clearly $ee'(e-e')=0$ and $(1-e)(1-e')(e-e')=0$. We have either $ee'\in A\setminus\mathfrak{p}$ or $(1-e)(1-e')\in A\setminus\mathfrak{p}$. It follows that $e-e'\in\Ker\pi_{\mathfrak{p}}$. Let $e-e'\in\mathfrak{J}$. Then $e-e'\in\bigcap
\limits_{\mathfrak{m}\in\Max(R)}\Ker\pi_{\mathfrak{m}}=0$. Here is another proof, we have $D(e)=D(e')$. Thus $e=e'$. $\Box$ \\

Using Lemma \ref{Lemma II}, then Corollary \ref{Corollary I} can also be phrased as: If $I$ is an ideal of a ring $R$ contained in the nilradical, then the ring map $\mathcal{B}(R)\rightarrow\mathcal{B}(R/I)$ given by $e\rightsquigarrow e+I$ is an isomorphism.

\begin{proposition}\label{Proposition II} Let $I\subseteq J$ be ideals of a ring $R$ such that the canonical map $\Max(R/J)\rightarrow\Max(R/I)$ is surjective. If $J$ is a lifting ideal, then $I$ is as well.
\end{proposition}

{\bf Proof.} If $f-f^{2}\in I$ for some $f\in R$, then there is an idempotent $e\in R$ such that $f-e\in J$. To prove $f-e\in I$, by Theorem \ref{Theorem II} it will be enough to show that $D(f+I)=D(e+I)$. Let $\mathfrak{p}$ be a prime ideal of $R$ containing $I$. There exists a maximal ideal $\mathfrak{m}$ of $R$ such that $\mathfrak{p}\subseteq\mathfrak{m}$. If $e\notin\mathfrak{p}$ then $1-e\in\mathfrak{p}$. If $f\in\mathfrak{p}$ then $e\in\mathfrak{m}$, since by hypothesis, $J\subseteq\mathfrak{m}$. But this is a contradiction. Thus $f\notin\mathfrak{p}$. To see the reverse inclusion, if $f\notin\mathfrak{p}$ then $1-f\in\mathfrak{p}$. If $e\in\mathfrak{p}$ then $f\in\mathfrak{m}$, because $J\subseteq\mathfrak{m}$. Again this is a contradiction. So $e\notin\mathfrak{p}$. $\Box$ \\

The converse of Proposition \ref{Proposition II} does not hold, see Example \ref{Example iii non lift}.

\begin{corollary} Let $I$ be an ideal of a ring $R$ such that $I\subseteq\mathfrak{J}$. If $\mathfrak{J}$ is a lifting ideal, then $I$ is as well.
\end{corollary}

{\bf Proof.} It follows from Proposition \ref{Proposition II}. $\Box$

\begin{corollary}\cite[Proposition 1.15]{Rostami et al.} Let $I$ be an ideal of a ring $R$. If $J=\bigcap\limits_{\mathfrak{m}\in V(I)\cap\Max(R)}\mathfrak{m}$ is a lifting ideal, then $I$ is as well.
\end{corollary}

{\bf Proof.} It follows from Proposition \ref{Proposition II}. $\Box$ \\

The following result is probably well known, but we provide a proof for the sake of completeness.

\begin{lemma}\label{Corollary III trivial idemps} If $R$ is either a local ring or has a unique minimal prime, then it has no nontrivial idempotents.
\end{lemma}

{\bf Proof.} If $R$ is a local ring the the assertion is clear. Assume $R$ has a unique minimal prime $\mathfrak{p}$. If $e\in R$ is an idempotent, then either $e\in\mathfrak{p}$ or $1-e\in\mathfrak{p}$. But $\mathfrak{p}$ is the nilradical of $R$. So $e\in\{0,1\}$. As another proof,
 $\Spec(R)=V(\mathfrak{p})$ is connected, since every irreducible space is connected. So $R$ has no nontrivial idempotents.  $\Box$

\begin{example}\label{Example iv idempotents} The converse of Lemma \ref{Corollary III trivial idemps} does not hold. For instance, we prove that $S=\mathbb{Z}[x]/I$ is an example of this type where $I=(x^{2}-1)$. First note that if $r_{1},\ldots,r_{n}$ are finite number of elements of an integral domain $R$, then $(x_{1}-r_{1},\ldots,x_{n}-r_{n})$ is a prime ideal of the polynomial ring $R[x_{1},\ldots,x_{n}]$, for its proof see e.g. \cite[Theorem 3.2]{Tarizadeh 3}. So $\mathfrak{p}/I$ and $\mathfrak{q}/I$ are the distinct minimal primes of $S$ where $\mathfrak{p}=(x+1)$ and $\mathfrak{q}=(x-1)$. Then we show that $S$ has no nontrivial idempotents. First note that if $R$ is a ring and $a\in R$, then we have the canonical isomorphism of rings $R[x]/(a,x+1)\simeq R/aR$. Now if $f-f^{2}\in I$ for some $f\in\mathbb{Z}[x]$, then either $f\in\mathfrak{p}$ or $1-f\in\mathfrak{p}$. If $f\in\mathfrak{p}$ then $f=(x+1)g\in\mathfrak{q}$, since $\mathfrak{p}+\mathfrak{q}=(2,x+1)\neq\mathbb{Z}[x]$. So $g\in\mathfrak{q}$. This yields that $f\in I$. If $1-f\in\mathfrak{p}$ then by an argument exactly like the above, we get that $1-f\in I$. Hence, $S$ has no nontrivial idempotents. Finally, if $p$ and $q$ are distinct prime numbers then $(p,x+1)$ and $(q,x+1)$ are distinct maximal ideals of $\mathbb{Z}[x]$. By passing to the quotient modulo $I$,
we find infinitely many maximal ideals in $S$. Hence, $S$ is not a local ring.
\end{example}

\begin{example}\label{Example iii non lift} The Jacobson radical of an integral domain with finitely many maximal ideals $\geqslant2$ is not a lifting ideal. In fact, let $R$ be a nonlocal ring with finitely many maximal ideals and a unique minimal prime ideal. Then the Jacobson radical of $R$ is not a lifting ideal. Because $R$ has at least two distinct maximal ideals, choose a maximal ideal $\mathfrak{m}$ of $R$ then $\mathfrak{m}+I=R$ where $I=\bigcap\limits_{\substack{\mathfrak{m}'\in\Max(R),\\
\mathfrak{m}'\neq\mathfrak{m}}}\mathfrak{m}'$. Thus there are $f\in\mathfrak{m}$ and $g\in I$ such that $f+g=1$. It follows that $f(1-f)=fg\in\mathfrak{J}$. If $\mathfrak{J}$ is a lifting ideal, then there is an idempotent $e\in R$ such that $f-e\in\mathfrak{J}$. By Lemma \ref{Corollary III trivial idemps}, $e\in\{0,1\}$. But this is impossible, since neither $f$ nor $1-f$ is in $\mathfrak{J}$. \\
\end{example}

Contrary to the nilradical which is always a lifting ideal (cf. Corollary \ref{Corollary I}), the Jacobson radical is not necessarily a lifting ideal, see Example \ref{Example iii non lift}. Investigation of the lifting property of this ideal leads us to the following general conclusion.

\begin{theorem}\label{Corollary iv nice} Let $R$ be a ring. Then the Zariski clopens of $\Max(R)$ are precisely of the form $\Max(R)\cap D(f)$ where $f(1-f)\in\mathfrak{J}$.
\end{theorem}

{\bf Proof.} Let $X=\Max(R)$. If $f(1-f)\in\mathfrak{J}$ then $X\cap D(f)=X\cap V(1-f)$ and so $X\cap D(f)$ is a Zariski clopen of $X$. Conversely, let $U$ be a Zariski clopen of $X$. It is well known that $X$ is quasi-compact. Thus we get that $U=\bigcup\limits_{i=1}^{n}X\cap D(f_{i})$, since every closed subset of a quasi-compact space is quasi-compact. Similarly, we may write $U^{c}=X\setminus U=\bigcup\limits_{k=1}^{m}X\cap D(g_{k})$. It follows that $X\cap D(f_{i})\cap D(g_{k})=\emptyset$ and so $f_{i}g_{k}\in\mathfrak{J}$ for all $i$ and $k$. Setting  $I=(f_{1},\ldots,f_{n})$ and $J=(g_{1},\ldots,g_{m})$. Then clearly $IJ\subseteq\mathfrak{J}$ and
$I+J=R$. Hence, $f+g=1$ for some $f\in I$ and $g\in J$. So $f(1-f)=fg\in IJ\subseteq\mathfrak{J}$. Clearly $X\cap D(f)\subseteq U$. If $\mathfrak{m}\in U$ then $f_{i}\notin\mathfrak{m}$ for some $i$. But
$f_{i}g_{k}\in\mathfrak{J}\subseteq\mathfrak{m}$ and so $g_{k}\in\mathfrak{m}$ for all $k$. This yields that $g\in\mathfrak{m}$ and so $f=1-g\notin\mathfrak{m}$. Hence, $\mathfrak{m}\in X\cap D(f)$.
Therefore $U=X\cap D(f)$. This completes the proof. $\Box$ \\

In the following result, $\Clop(X)$ denotes the Zariski clopens of $X=\Max(R)$. This result is the analogue of \cite[Theorem 3.1]{Tarizadeh 2}.

\begin{corollary}\label{Theorem v new} Let $R$ be a ring and $X=\Max(R)$. Then the map $\phi:\mathcal{B}(R/\mathfrak{J})\rightarrow\Clop(X)$ given by $f+\mathfrak{J}\rightsquigarrow X\cap D(f)$ is an isomorphism of rings.
\end{corollary}

{\bf Proof.} It is a morphism of rings, because $\phi=\Clop(\eta)\circ\psi$ where the canonical continuous map $\eta:X\rightarrow Y=\Spec(R/\mathfrak{J})$ is given by $\mathfrak{m}\rightsquigarrow\mathfrak{m}/\mathfrak{J}$ and for the canonical isomorphism of rings $\psi:\mathcal{B}(R/\mathfrak{J})\rightarrow
\Clop(Y)$ see \cite[Theorem 3.1]{Tarizadeh 2}. This map is also injective, since if $X\cap D(f)=\emptyset$ then $f\in\mathfrak{J}$. Finally by Theorem \ref{Corollary iv nice}, it is surjective. $\Box$

\begin{corollary}\label{Corollary v Jacobson} Let $R$ be a ring and  $X=\Max(R)$. Then the Jacobson radical of $R$ is a lifting ideal if and only if the ring map $\mathcal{B}(R)\rightarrow\Clop(X)$ given by $e\rightsquigarrow X\cap D(e)$ is an isomorphism.
\end{corollary}

{\bf Proof.} The ring map $\mathcal{B}(R)\rightarrow\mathcal{B}(R/\mathfrak{J})$ given by $e\rightsquigarrow e+\mathfrak{J}$ is always injective. Thus by Corollary \ref{Theorem v new}, the ring map $\mathcal{B}(R)\rightarrow\Clop(X)$ given by $e\rightsquigarrow X\cap D(e)$ is also injective. Now under the light of Theorem \ref{Corollary iv nice}, this ring map is surjective
(and hence isomorphism) if and only if the Jacobson radical of $R$ is a lifting ideal. $\Box$ \\

If $I$ is an ideal of a ring $R$, then the ideal $(f\in I: f=f^{2})$ is denoted by $I^{\ast}$. This ideal is quite interesting, for more information and applications see \cite[Theorem 5.1(xii)]{Tarizadeh}, \cite{Magid}, \cite{Rostami et al.}, \cite[\S4]{Tarizadeh 2} and \cite[Proposition 3.11]{Tarizadeh 5}. The following result simplifies
and improves on \cite[Theorems 1.5 and 3.2]{Rostami et al.}, with the addition of conditions (iv) and (v).

\begin{theorem}\label{Theorem vi new} Let $I$ be an ideal of a ring $R$. Then the following statements are equivalent. \\
$\mathbf{(i)}$ $I$ is a lifting ideal of $R$. \\
$\mathbf{(ii)}$ If $J_{1}J_{2}\subseteq I$ for some coprime ideals
$J_{1}$ and $J_{2}$ of $R$, then there exists an idempotent $e\in R$ such that $e\in I+J_{1}$ and $1-e\in I+J_{2}$. \\
$\mathbf{(iii)}$ If $I=J_{1}J_{2}$ for some coprime ideals
$J_{1}$ and $J_{2}$ of $R$, then there exists an idempotent $e\in R$ such that $e\in J_{1}$ and $1-e\in J_{2}$. \\
$\mathbf{(iv)}$ The canonical ring map $R\rightarrow R[x]/J$ lifts idempotents where $J$ is the kernel of the ring map $R[x]\rightarrow R/I$ which is given by $f(x)\rightsquigarrow f(0)+I$. \\
$\mathbf{(v)}$ If $\mathfrak{m}$ is a maximal ideal of $R$, then the quotient ring $R/(I+\mathfrak{m}^{\ast})$ has no nontrivial idempotents.
\end{theorem}

{\bf Proof.} $\mathbf{(i)}\Rightarrow\mathbf{(ii)}:$ There are $f\in J_{1}$ and $g\in J_{2}$ such that $f+g=1$. This yields that $f-f^{2}=fg\in I$. Thus by the hypothesis, there is an idempotent $e\in R$ such that $e-f\in I$. \\
$\mathbf{(ii)}\Rightarrow\mathbf{(iii)}:$ There is nothing to prove. \\
$\mathbf{(iii)}\Rightarrow\mathbf{(i)}:$ Assume $f-f^{2}\in I$ for some $f\in R$. Clearly $J_{1}=Rf+I$ and $J_{2}=R(1-f)+I$ are coprime ideals. Also $I=J_{1}J_{2}$. Because clearly $J_{1}J_{2}\subseteq I$. Conversely, if $g\in I$ then we may write $g=g(1-f)+fg\in J_{1}J_{2}$. Thus by hypothesis, there is an idempotent $e\in R$ such that $e\in J_{1}$ and $1-e\in J_{2}$. Then clearly $e-f\in J_{1}$ and also $e-f=(e-1)+(1-f)\in J_{2}$. Thus $e-f\in J_{1}\cap J_{2}=J_{1}J_{2}=I$. \\
$\mathbf{(i)}\Rightarrow\mathbf{(iv)}:$ If $f-f^{2}\in J$ for some $f=\sum\limits_{i=0}^{n}r_{i}x^{i}\in R[x]$, then $r_{0}-r^{2}_{0}\in I$. Thus by the hypothesis, there is an idempotent $e\in R$ such that $r_{0}-e\in I$. It follows that $f-e\in J$. \\
$\mathbf{(iv)}\Rightarrow\mathbf{(i)}:$ Clear. \\
$\mathbf{(i)}\Rightarrow\mathbf{(v)}:$ If $f-f^{2}\in I+\mathfrak{m}^{\ast}$ for some $f\in R$, then there are $g\in I$ and an idempotent $e\in\mathfrak{m}$ such that $f-f^{2}=g+re$ for some $r\in R$. This yields that $f(1-e)-\big(f(1-e)\big)^{2}=(f-f^{2})(1-e)=(1-e)g\in I$. Thus be hypothesis, there is an idempotent $e'\in R$ such that $f(1-e)-e'\in I$. So $f=h+fe+e'$ for some $h\in I$.
We have either $e'\in\mathfrak{m}$ or $1-e'\in\mathfrak{m}$. If $e'\in\mathfrak{m}$ then $fe+e'\in\mathfrak{m}^{\ast}$ and so $f\in I+\mathfrak{m}^{\ast}$. If $1-e'\in\mathfrak{m}$ then $1-e'-fe\in\mathfrak{m}^{\ast}$ and so
$1-f=1-e'-fe-h\in I+\mathfrak{m}^{\ast}$. \\
$\mathbf{(v)}\Rightarrow\mathbf{(i)}:$ If $f(1-f)\in I$ for some $f\in R$, then $(I:f)^{\ast}+\big(I:(1-f)\big)^{\ast}$ is the whole ring $R$. If not, then it is contained in a maximal ideal $\mathfrak{m}$. Then by hypothesis, either $f\in I+\mathfrak{m}^{\ast}$ or $1-f\in I+\mathfrak{m}^{\ast}$. If $f\in I+\mathfrak{m}^{\ast}$ then there exist some $g\in I$ and an idempotent $e\in\mathfrak{m}$ such that $f=g+re$ for some $r\in R$. It follows that $1-e\in (I:f)^{\ast}\subseteq\mathfrak{m}$ which is a contradiction. If $1-f\in I+\mathfrak{m}^{\ast}$ then we reach a contradiction by a similar argument. Therefore $(I:f)^{\ast}+\big(I:(1-f)\big)^{\ast}=R$. Thus there are idempotents $e_{1}\in(I:f)$ and $e_{2}\in\big(I:(1-f)\big)$ such that $1=r_{1}e_{1}+r_{2}e_{2}$ for some $r_{1},r_{2}\in R$. It follows that $f(1-e_{2})=r_{1}fe_{1}(1-e_{2})\in I$ and so
$f-e_{2}=f(1-e_{2})-e_{2}(1-f)\in I$. $\Box$ \\

If $I+J$ is a lifting ideal of a ring $R$, then the summands $I$ and $J$ are not necessarily lifting ideals. For example, $6\mathbb{Z}+35\mathbb{Z}=\mathbb{Z}$ is a lifting ideal of $\mathbb{Z}$, but none of the summands is a lifting ideal.
Also, if $I$ and $J$ are lifting ideals of a ring $R$, then $I+J$ is not necessarily a lifting ideal, see e.g. \cite[\S3]{Diesl et al.}. In this direction, we obtain Corollary \ref{Corollary vi regular} that generalizes a number of related results in the literature, especially \cite[Lemma 1.11(3)]{Rostami et al.} which states that every regular ideal is a lifting ideal. Indeed, let $I$ be a regular ideal of a ring $R$. If $f-f^{2}\in I$ for some $f\in R$, then $f-f^{2}=re$ for some idempotent $e\in I$. It follows that $(f-f^{2})(1-e)=0$. So $f(1-e)=f^{2}(1-e)$. Hence, $f(1-e)$ is
an idempotent and $f-f(1-e)=fe\in I$. More generally we have the following result.

\begin{corollary}\label{Corollary vi regular} If $I$ is a lifting ideal of a ring $R$ and $J$ is a regular ideal of $R$, then $I+J$ is a lifting ideal of $R$.
\end{corollary}

{\bf Proof.} By Theorem \ref{Theorem vi new}(v), it will be enough to show that for each $\mathfrak{m}\in\Max(R)$ then $R/(I+J+\mathfrak{m}^{\ast})$ has no nontrivial idempotents. If $f(1-f)\in I+J+\mathfrak{m}^{\ast}$ for some $f\in R$, then there exist $g\in I$ and an idempotent $e\in J+\mathfrak{m}^{\ast}$ such that $f(1-f)=g+re$ for some $r\in R$. It follows that $f(1-e)-\big(f(1-e)\big)^{2}=g(1-e)\in I$. Thus by hypothesis, there exists an idempotent $e'\in R$ such that $f(1-e)-e'\in I$. So $f=h+fe+e'$ for some $h\in I$.
We have either $e'\in\mathfrak{m}$ or $1-e'\in\mathfrak{m}$. If $e'\in\mathfrak{m}$ then $fe+e'\in J+\mathfrak{m}^{\ast}$ and so $f\in I+J+\mathfrak{m}^{\ast}$. If $1-e'\in\mathfrak{m}$ then we get that $1-f\in I+J+\mathfrak{m}^{\ast}$. $\Box$ \\

Remember that an ideal $I$ of a ring $R$ is called a pure ideal if whenever $f\in I$ then there exists some $g\in I$ such that $f(1-g)=0$.
Every regular ideal is both pure and lifting ideal. But the converse does not hold. Finding a specific example of a pure (and lifting) ideal which is not regular ideal is not easy at all, but fortunately \cite[Example 4.7]{Tarizadeh 4} provides an implicit example of this type. Note that both a prime ideal $\mathfrak{p}$ and $\Ker\pi_{\mathfrak{p}}$ are lifting ideals of a ring $R$, since $R/\mathfrak{p}$ and $R/\Ker\pi_{\mathfrak{p}}$ have no nontrivial idempotents.

\begin{remark}\label{Remark I} Remember that a ring $R$ is called a clean ring if each $f\in R$ can be written as $f=e+g$ where $e\in R$ is an idempotent and $g\in R$ is invertible in $R$. This simple notion has several nontrivial equivalents, see e.g. \cite[\S5]{Tarizadeh}. Especially \cite[Theorem 5.1(ix)]{Tarizadeh} states that a ring $R$ is a clean ring if and only if for any distinct maximal ideals $\mathfrak{m}$ and $\mathfrak{m}'$ of $R$ there exists an idempotent $e\in R$ which lies in $\mathfrak{m}$ but not in $\mathfrak{m}'$, (note that the proof of the implication ``$\Rightarrow$" is easy, but the reverse implication is slightly technical). Using this, in Proposition \ref{Theorem III} we give a new and simple proof to the main result of \cite{Nicholson} in the commutative case.
\end{remark}

\begin{proposition}\label{Theorem III} Let $R$ be a ring. Then $R$ is a clean ring if and only if each ideal of $R$ is a lifting ideal.
\end{proposition}

{\bf Proof.} Assume $I$ is an ideal of a clean ring $R$ and $f\in R$ such that $f-f^{2}\in I$. We may write $f=e+g$ where $e\in R$ is an idempotent and $g\in R$ is invertible in $R$. Clearly $e'-f=g^{-1}(f-f^{2})\in I$ where $e'=1-e$. Conversely, if $\mathfrak{m}$ and $\mathfrak{m}'$ are distinct maximal ideals of $R$ then there exist  $f\in\mathfrak{m}$ and $g\in\mathfrak{m}'$ such that $f+g=1$. It follows that $f-f^{2}=fg\in\mathfrak{m}\mathfrak{m}'$. So by the hypothesis, there exists an idempotent $e\in R$ such that $f-e\in\mathfrak{m}\mathfrak{m}'$. This yields that $e\in\mathfrak{m}\setminus\mathfrak{m}'$. $\Box$ \\

There is further fundamental result on lifting idempotents which is the dual of Proposition \ref{Theorem III}. It states that if the idempotents of a ring $R$ can be lifted along each localization, then $R$ is a purified ring (i.e., for any distinct minimal primes $\mathfrak{p}$ and $\mathfrak{q}$ of $R$, there exists an idempotent which lies in $\mathfrak{p}$ but not in $\mathfrak{q}$). If moreover $R$ is a reduced ring, then the converse holds. For more information see \cite[Theorem 8.5]{Tarizadeh}. \\

If $\mathfrak{p}$ and $\mathfrak{q}$ are distinct minimal primes of a ring $R$, then there exist $f\in R\setminus\mathfrak{p}$ and $g\in R\setminus\mathfrak{q}$ such that $fg=0$ (see e.g. \cite[Lemma 3.2]{Tarizadeh}). It follows that the induced Zariski topology over the minimal spectrum $\Min(R)$ is Hausdorff. In particular, a ring $R$ is zero dimensional if and only if $\Spec(R)$ is Hausdorff. The following result gives us new characterizations of zero dimensional rings which are needed in the sequel. For further criteria see Remark \ref{Remark iii reverse}, \cite[Theorem 3.3]{Tarizadeh} and \cite[Theorem 2.2]{Tarizadeh 4}.

\begin{lemma}\label{Lemma vi zero dim} For a ring $R$ the following statements are equivalent. \\
$\mathbf{(i)}$ $\dim(R)=0$. \\
$\mathbf{(ii)}$ If $\mathfrak{p}$ and $\mathfrak{q}$ are distinct primes of $R$, then there exists an idempotent $e\in\mathfrak{p}\setminus\mathfrak{q}$. \\
$\mathbf{(iii)}$ $\Spec(R)$ is totally disconnected. \\
$\mathbf{(iv)}$ $\Spec(R)$ is homeomorphic to $\Spec\big(\mathcal{B}(R)\big)$.
\end{lemma}

{\bf Proof.} $\mathbf{(i)}\Rightarrow\mathbf{(ii)}:$ We may choose some $f\in\mathfrak{p}$ such that $f\notin\mathfrak{q}$. It is well known that $R/\mathfrak{N}$ is an absolutely flat (von Neumann regular) ring. Thus $f(1-fh)\in\mathfrak{N}$ for some $h\in R$. So there exists a natural number $n\geqslant1$ such that $f^{n}(1-fg')=0$ for some $g'\in R$. This yields that $f^{n}=f^{2n}g$ for some $g\in R$. So $e:=f^{n}g$ is an idempotent and $e\in\mathfrak{p}\setminus\mathfrak{q}$. \\
$\mathbf{(ii)}\Rightarrow\mathbf{(iii)}:$ Assume $\mathfrak{p}$ and $\mathfrak{q}$ are distinct points of a connected subset $C$ of $\Spec(R)$. By hypothesis, there exists an idempotent $e\in\mathfrak{p}\setminus\mathfrak{q}$. Clearly $\mathfrak{p}\in U=D(1-e)\cap C$ and $\mathfrak{q}\in V=D(e)\cap C$. We have $C=U\cup V$, since $e$ is an idempotent. But this is in contradiction with connectedness of $C$. Hence, $C$ is singleton. \\
$\mathbf{(iii)}\Rightarrow\mathbf{(iv)}:$ By hypothesis and \cite[Theorem 4.1]{Tarizadeh 2}, $\Spec\big(\mathcal{B}(R)\big)\simeq\Sp(R)\simeq\pi_{0}
\big(\Spec(R)\big)\simeq\Spec(R)$. \\
$\mathbf{(iv)}\Rightarrow\mathbf{(i)}:$ If $\mathfrak{p}\subseteq\mathfrak{q}$ are primes of $R$, then $\mathfrak{q}\in V(\mathfrak{p})=\overline{\{\mathfrak{p}\}}$. By hypothesis, there exists a homeomorphism
$\phi:\Spec(R)\rightarrow\Spec\big(\mathcal{B}(R)\big)$. It follows that $\phi(\mathfrak{q})\in\overline{\{\phi(\mathfrak{p})\}}=
V\big(\phi(\mathfrak{p})\big)=\{\phi(\mathfrak{p})\}$, since in a Boolean ring each prime is maximal. Therefore $\phi(\mathfrak{p})=\phi(\mathfrak{q})$ and so $\mathfrak{p}=\mathfrak{q}$. $\Box$ \\

The following result follows from Proposition \ref{Theorem III} and \cite[Theorem 8.5(vi)]{Tarizadeh}. But here we provide a new proof which uses less machinery.

\begin{corollary}\label{Corollary viii along locz} If $R$ is a zero dimensional ring, then the idempotents of $R$ can be lifted modulo each ideal. If moreover $R$ is reduced, then the idempotents of $R$ can be lifted along each localization.
\end{corollary}

{\bf Proof.} Let $I$ be an ideal of $R$ and $f-f^{2}\in I$ for some $f\in R$. There exists a natural number $n\geqslant1$ such that $f^{n}=f^{2n}g$ for some $g\in R$, see the proof of Lemma \ref{Lemma vi zero dim}. Thus $f^{n}g$ is an idempotent. Clearly $D(f+I)=D(f^{n}g+I)$, in fact $D(f)=D(f^{n}g)$. So by Theorem \ref{Theorem II}, $f-f^{n}g\in I$. Now assume $R$ is reduced. If $f/s\in S^{-1}R$ is an idempotent for some multiplicative subset $S$ of $R$, then $s'(fs^{2}-f^{2}s)=0$ for some $s'\in S$. We also have $f=f^{2}g$ for some $g\in R$, since $R$ is absolutely flat. Thus $fg$ is an idempotent and $ss'(f-fgs)=0$, so $fg/1=f/s$. $\Box$

\begin{corollary}\cite[Proposition 1.5]{Nicholson} Let $I$ be a lifting ideal of a ring $R$ contained in the Jacobson radical. If $R/I$ is a clean ring, then $R$ is so.
\end{corollary}

{\bf Proof.} It follows from Remark \ref{Remark I}. $\Box$

\begin{lemma}\label{Lemma I} Let $I$ be an ideal of a ring $R$ and $f\in R$ such that $f-f^{2}\in I^{k}$ for some $k\geqslant1$. If $n\geqslant k$ then there exists $g\in R$ such that $g-g^{2}\in I^{n}$ and $f-g\in I^{k}$.
\end{lemma}

{\bf Proof.} By Corollary 3.7, the canonical ring map $R/I^{n}\rightarrow R/I^{k}$ lifts idempotents. $\Box$ \\

We provide new and alternative proofs to the following results.

\begin{proposition}\label{Theorem viii complete}\cite[Theorem 21.31]{Lam} Let $I$ be an ideal of a ring $R$. If $R$ is complete with respect to the $I$-adic topology, then $I$ is a lifting ideal.
\end{proposition}

{\bf Proof.} Assume $f_{1}-f^{2}_{1}\in I$ for some $f_{1}\in R$. By Lemma \ref{Lemma I}, there exists $f_{2}\in R$ such that $f_{2}-f^{2}_{2}\in I^{2}$ and $f_{2}-f_{1}\in I$. In fact, by the successive applications of Lemma \ref{Lemma I}, then for each $n\geqslant2$ we may find $f_{n}\in R$ such that $f_{n}-f^{2}_{n}\in I^{n}$ and $f_{n}-f_{n-1}\in I^{n-1}$. It follows that the sequence $(f_{n}+I^{n})_{n\geqslant1}$ is an idempotent and it is a member of the $I$-adic completion $\widehat{R}=
\lim\limits_{\overleftarrow{n\geqslant1}}R/I^{n}$. Thus by the hypothesis, there is an idempotent $e\in R$ such that $f_{n}-e\in I^{n}$ for all $n\geqslant1$. $\Box$

\begin{proposition}\cite[Proposition 1.2]{Rostami et al.}\label{Prop RHetal} Let $I$ and $J$ be non-coprime ideals of a ring $R$ such that $I$ is a lifting ideal and $R/J$ has no nontrivial idempotents. Then $I\cap J$ is a lifting ideal.
\end{proposition}

{\bf Proof.} If $f-f^{2}\in I\cap J$ for some $f\in R$, then there is an idempotent $e\in R$ such that $e-f\in I$. It suffices to show that $e-f\in J$. There exists a prime ideal $\mathfrak{p}$ of $R$ such that $I+J\subseteq\mathfrak{p}$. We have either $f\in J$ or $1-f\in J$. If $f\in J$ then $e\in\mathfrak{p}$ and so $e\in J$. Thus $e-f\in J$. If $1-f\in J$ then $1-e\in\mathfrak{p}$ and so $1-e\in J$. Hence, $e-f=(e-1)+(1-f)\in J$. $\Box$

\section{Bounding a Boolean ring by its primes}

Theorems \ref{Theorem Deligne-Tarizadeh} and \ref{Corollary ix main} are the main results of this section. In this section, in particular, we investigate how many idempotents can a commutative ring have?

\begin{lemma}\label{Lemma v infinite cardinal} For any infinite cardinal $\kappa$, there is a Boolean ring $R$ with $|R|=|\Spec(R)|=\kappa$.
\end{lemma}

{\bf Proof.} Let $X$ be an infinite set and let $R$ be the Boolean ring of all subsets of $X$ which are either finite or cofinite (i.e., its complement is finite). In fact, $R$ is a subring of the power set ring $\mathcal{P}(X)$.
One can check that the prime (maximal) ideals of $R$ are precisely either $\Fin(X)$ or of the form $\mathfrak{m}_{x}\cap R$ where $x\in X$ and $\mathfrak{m}_{x}=\mathcal{P}(X\setminus\{x\})$. So $|\Spec(R)|=|X|+1=|X|$, since $X$ is infinite. It is well known that $|X|=|\Fin(X)|$. So $|X|\leqslant |R|$, since $\Fin(X)\subset R$. Finally, consider the map $\phi:R\rightarrow\Fin(X)\times\{0,1\}$ where $\phi(A)$ is either $(A,0)$ or $(A^{c},1)$, according as $A$ or $A^{c}=X\setminus A$ is finite. Clearly $\phi$ is injective, so $|R|\leqslant|\Fin(X)|$. Therefore $|R|=|X|=|\Spec(R)|$. $\Box$

\begin{proposition}\label{proposition vi zero dim} Let $R$ be a ring. If $|R|\leqslant|\Spec(R)|$, then $\Spec(R)$ is infinite.
\end{proposition}

{\bf Proof.} If $\Spec(R)$ is finite then $R$ is as well. Hence, $R$ is an Artinian ring. Therefore by the structure theorem for Artinian rings (cf. Remark \ref{Example II}),  $|R|\geqslant|\mathcal{B}(R)|=2^{n}>n=|\Spec(R)|$ which is a contradiction. $\Box$ \\

The converse of the above result does not hold. For example, take the ring $B\times F$ where $B$ is a countably infinite Boolean ring with $|B|=|\Spec(B)|$, (cf. Lemma \ref{Lemma v infinite cardinal}) and $F$ is an uncountable field. But if $R$ is Boolean, then the converse of Proposition \ref{proposition vi zero dim} holds, see Theorem \ref{Theorem Deligne-Tarizadeh}.

\begin{remark} Let $R$ be a ring such that $|R|<|\Spec(R)|$. Then by Proposition \ref{proposition vi zero dim}, $R$ is infinite. Clearly $|\Spec(R)|\leqslant2^{|R|}$, because $\Spec(R)$ is a subset of $\mathcal{P}(R)$. Therefore by the Generalized Continuum Hypothesis,
 $|\Spec(R)|=2^{|R|}$.
\end{remark}

Let $R$ be a nonzero ring. Then the polynomial ring $R[x_{1},\ldots,x_{n}]$ with finitely many variables has the cardinality either $\aleph_{0}$ or $|R|$, according as $R$ is finite or infinite. In a formula, $R[x_{1},\ldots,x_{n}]$ has the cardinality $|R|\aleph_{0}$.
It also holds for the polynomial ring $R[x_{1},x_{2},x_{3},\ldots]$ with countably infinite variables.

\begin{lemma}\label{Lemma I Tari} Let $R$ be a nonzero ring and $I$ an infinite set. Then the cardinality of the polynomial ring $R[x_{i}: i\in I]$ is $|R\times I|$.
\end{lemma}

{\bf Proof.} The map $R\times I\rightarrow S=R[x_{i}: i\in I]$ given by $(r,i)\rightsquigarrow r+x_{i}$ is injective. Hence, $S$ has the cardinality $\geqslant|R\times I|$. To complete the proof we act as follows. It is well known that $|\Fin(I)|=|I|$. The set $\Fin(I)$ ordered by the inclusion is directed.
If $J\in\Fin(I)$ then the polynomial ring $R[x_{i}: i\in J]$ with finitely many variables has cardinality $|R|\aleph_{0}$.
The polynomial rings $R[x_{i}: i\in J]$ with $J\in\Fin(I)$ together with the canonical injections, as the transition morphisms, is a direct system of rings over the poset $\Fin(I)$. The ring $S$ is canonically isomorphic to the direct limit of this system. Now $S$ being a quotient of a disjoint union of the $R[x_{i}: i\in J]$ indexed by $\Fin(I)$, it has cardinality
$\leqslant|R|\aleph_{0}|I|=|R|.|I|=|R\times I|$. $\Box$ \\

We obtained the following result independently of \cite[Chap. II, \S5, Theorem 11, Exercise 36]{Gratzer} and \cite[Lemma 11.1]{R.S. Pierce}.

\begin{theorem}\label{Theorem Deligne-Tarizadeh} If $R$ is an infinite  Boolean ring, then $|R|\leqslant|\Spec(R)|$.
\end{theorem}

{\bf Proof.} Let $X=\Spec(R)$ and  $\Delta=\{(\mathfrak{p},\mathfrak{p}):\mathfrak{p}\in X\}$ be the diagonal. If $(\mathfrak{p},\mathfrak{q})\in X^{\ast}:=X^{2}\setminus\Delta$, then there is some $e_{\mathfrak{p},\mathfrak{q}}\in R$ such that $e_{\mathfrak{p},\mathfrak{q}}\in\mathfrak{p}\setminus\mathfrak{q}$.
We claim that $R$ as $\mathbb{Z}_{2}-$algebra is generated by the $e_{\mathfrak{p},\mathfrak{q}}$'s. That is,
$R=\mathbb{Z}_{2}[e_{\mathfrak{p},\mathfrak{q}}: (\mathfrak{p},\mathfrak{q})\in X^{\ast}]$. Clearly $S=\mathbb{Z}_{2}[e_{\mathfrak{p},\mathfrak{q}}: (\mathfrak{p},\mathfrak{q})\in X^{\ast}]$ is a subring of $R$, since $\mathbb{Z}_{2}=\{0,1\}$ is a subring of $R$. Take $e\in R$. If $\mathfrak{p}\in D(e)$ then $V(e)\subseteq\bigcup\limits_{\mathfrak{q}\in V(e)}D(e_{\mathfrak{p},\mathfrak{q}})$. Using the compactness of $V(e)$, then we may find a finite subset $\{\mathfrak{q}_{i}:i=1,\ldots,n_{\mathfrak{p}}\}\subseteq V(e)$ such that $V(e)\subseteq\bigcup\limits_{i=1}^{n_{\mathfrak{p}}}
D(e_{\mathfrak{p},\mathfrak{q}_{i}})$. If follows that $\mathfrak{p}\in D(s_{\mathfrak{p}})\subseteq D(e)$ where $s_{\mathfrak{p}}=\prod\limits_{i=1}^{n_{\mathfrak{p}}}
(1-e_{\mathfrak{p},\mathfrak{q}_{i}})\in S$. Now using the compactness of $D(e)$, then we may find a finite subset $\{\mathfrak{p}_{k}:i=1,\ldots,d\}\subseteq D(e)$ such that $D(e)=\bigcup\limits_{k=1}^{d}D(s_{\mathfrak{p}_{k}})$. Using Theorem \ref{Theorem II}, then $e=1-\prod\limits_{k=1}^{d}(1-s_{\mathfrak{p}_{k}})\in S$. This establishes the claim. Therefore $R$ is a quotient of the polynomial ring $T=\mathbb{Z}_{2}[x_{\mathfrak{p},\mathfrak{q}}:
(\mathfrak{p},\mathfrak{q})\in X^{\ast}]$. By the hypothesis, $R$ is infinite and so by Theorem \ref{Theorem II}, $X$ is also infinite. But for any infinite cardinal $\kappa$, then $\kappa.\kappa=\kappa$ (see e.g. \cite[p. 162, Lemma 6R]{Enderton}). It follows that $|X^{\ast}|=|X|$.
Now using Lemma \ref{Lemma I Tari}, then $|R|\leqslant|T|=|X^{\ast}|=|X|$. $\Box$

\begin{example} There are Boolean rings such that the inequality of Theorem \ref{Theorem Deligne-Tarizadeh} is strict for them.
As an example, take the power set ring $\mathcal{P}(X)$ where $X$ is an infinite set with the cardinality $\kappa$. Remember that $\Spec\big(\mathcal{P}(X)\big)$ is the Stone-\v{C}ech compactification of the discrete space $X$. So $|\mathcal{P}(X)|=2^{\kappa} <|\Spec\big(\mathcal{P}(X)\big)|=2^{2^{\kappa}}$.
\end{example}

\begin{lemma}\label{Theorem ix finiteness} The following statements are equivalent for a ring $R$. \\
$\mathbf{(i)}$ $R$ has finitely many idempotents. \\
$\mathbf{(ii)}$ $\Spec(R)$ has finitely many Zariski connected components.\\
$\mathbf{(iii)}$ The Boolean ring $\mathcal{B}(R)$ has finitely many prime ideals. \\
In this case, $|\mathcal{B}(R)|=2^{n}$ where $n=|\pi_{0}\big(\Spec(R)\big)|=|\Spec\big(\mathcal{B}(R)\big)|$.
\end{lemma}

{\bf Proof.} Under the light of homeomorphisms $\Spec\big(\mathcal{B}(R)\big)\simeq\Sp(R)
\simeq\pi_{0}\big(\Spec(R)\big)$, then the implications $\mathbf{(i)}\Rightarrow\mathbf{(ii)}\Rightarrow\mathbf{(iii)}$ are clear. \\
$\mathbf{(iii)}\Rightarrow\mathbf{(i)}:$ By Theorem \ref{Theorem II}, $\mathcal{B}(R)$ is finite. \\
Now if one of these equivalent conditions, say (iii), holds then
$X=\Spec\big(\mathcal{B}(R)\big)$ is discrete, since every finite Hausdorff space is discrete. Therefore $\mathcal{B}(R)\simeq\Clop(X)=\mathcal{P}(X)$. So
$|\mathcal{B}(R)|=2^{n}$ where $n=|X|=|\pi_{0}\big(\Spec(R)\big)|$.
$\Box$

\begin{theorem}\label{Corollary ix main} Let $R$ be a ring and $\kappa=|\pi_{0}\big(\Spec(R)\big)|$. Then: $$|\mathcal{B}(R)|<\infty\Leftrightarrow|\mathcal{B}(R)|=2^{\kappa}.$$
\end{theorem}

{\bf Proof.} The implication ``$\Rightarrow$" follows from Lemma \ref{Theorem ix finiteness}. Conversely, if $\mathcal{B}(R)$ is infinite then by Theorem \ref{Theorem Deligne-Tarizadeh}, $2^{\kappa}=|\mathcal{B}(R)|\leqslant|\Spec\big(\mathcal{B}(R)\big)|=
\kappa$. But this is in contradiction with Cantor's theorem (i.e., $\kappa<2^{\kappa}$ for all cardinals $\kappa$). Hence, $\mathcal{B}(R)$ is finite. $\Box$ \\

In regarding with Theorem \ref{Corollary ix main}, note that if $|\mathcal{B}(R)|=2^{\kappa}$ for some cardinal $\kappa$, then the number of idempotents of $R$ is not necessarily finite. For example, take $\mathcal{P}(X)$ where $X$ is an infinite set.

\begin{corollary} Let $R$ be a zero dimensional ring and $\kappa=|\Spec(R)|$. Then $\mathcal{B}(R)$ is finite if and only if $|\mathcal{B}(R)|=2^{\kappa}$.
\end{corollary}

{\bf Proof.} It follows from Theorem \ref{Corollary ix main} and Lemma \ref{Lemma vi zero dim}(iii). $\Box$ \\

The following result can be viewed as a generalization of Lemma \ref{Corollary III trivial idemps}.

\begin{proposition}\label{Prop v mini-max} Let $R$ be a ring. If either $\Min(R)$ or $\Max(R)$ is finite,
then $R$ has finitely many idempotents.
\end{proposition}

{\bf Proof.} The ring $R/\mathfrak{N}$ can be viewed as a subring of $\prod\limits_{\mathfrak{p}\in\Min(R)}R/\mathfrak{p}$.
If $\Min(R)$ is finite, then $R/\mathfrak{N}$ is a subring of a finite product of integral domains. Hence, $R/\mathfrak{N}$ and so $R$ have finitely many idempotents. If $\Max(R)$ is finite, then by Chinese remainder theorem, $R/\mathfrak{J}$ is a finite product of fields. Thus $R/\mathfrak{J}$ has finitely many idempotents. But for any ring $R$, the ring map $\mathcal{B}(R)\rightarrow\mathcal{B}(R/\mathfrak{J})$ given by $e\rightsquigarrow e+\mathfrak{J}$ is injective. Hence, $R$ has finitely many idempotents. $\Box$

\begin{example} The converse of Proposition \ref{Prop v mini-max} does not hold. For example, consider $S=\mathbb{Z}[x]/(x^{2}-1)$ in Example \ref{Example iv idempotents}. Let $T$ be a prime-inverse ring of $S$ (for more information see \cite[Remark 3.4 and Corollary 3.5]{Tarizadeh 2}).
Then the ring $R:=S\times T$ has exactly 4 idempotents but neither its minimal primes nor its maximal ideals is finite.
\end{example}

\begin{remark} The identification  $\Spec\big(\mathcal{B}(R)\big)\simeq\Sp(R)$ is so useful to compute the connected components of some complicated prime spectra in a more explicit way. For instance, let $(R_{x})$ be a family of rings indexed by a set $X$ such that each $R_{x}$ has no nontrivial idempotents (each local ring or integral domain is a such ring). Let $R=\prod\limits_{x\in X}R_{x}$. Then the canonical map $\mathcal{B}(R)\rightarrow\mathcal{P}(X)$ given by $(e_{x})\rightsquigarrow\{x\in X: e_{x}=1\}$
is an isomorphism of rings. It follows that $\pi_{0}\big(\Spec(R)\big)$ is homeomorphic to $\Spec\big(\mathcal{P}(X)\big)$, the Stone-\v{C}ech compactification of the discrete space $X$. In particular, if $\kappa=|X|$, then the number of connected components of $\Spec(R)$ is either $\kappa$ or $2^{2^{\kappa}}$, according as $X$ is finite or infinite.
\end{remark}

\begin{remark} If $e$ is an idempotent of a ring $R$, then $(1-2e)^{2}=1$. So $1-2e$ is invertible in $R$. In fact, the map $\eta(e)=1-2e$ is a morphism of groups from the additive group $\mathcal{B}(R)$ to the group of units of $R$. That is, $\eta(e\oplus e')=(1-2e)(1-2e')$. Clearly $\Ker\eta=\{e\in\mathcal{B}(R): 2e=0\}$ is an ideal of the Boolean ring $\mathcal{B}(R)$. If $2$ is invertible in $R$, then $\eta$ is injective and its image is the subgroup $\{f \in R: f^{2}=1\}$, because if $f^{2}=1$ then $e:=(1-f)/2$ is an idempotent. But the converse does not hold. For example, take $R=\mathbb{Z}$ or $R=\mathbb{Z}/4\mathbb{Z}$. Moreover, the ideal $I=\Ker\eta$ is not necessarily a prime ideal of $\mathcal{B}(R)$. For example, if $R=\mathbb{Z}/30\mathbb{Z}$ then $I=\{0,15\}$ is not a prime ideal of $\mathcal{B}(R)=\{0,1,6,10,15,16,21,25\}$, since $6\times 10=0\in I$ but neither $6$ nor $10$ is in $I$.
\end{remark}

\section{Orthogonal and primitive idempotents}

Theorems \ref{Theorem  vii orthogonality}, \ref{Theorem lift h0},
\ref{Theorem x isolated points}, \ref{theorem xi primitive} and \ref{Theorem xii TD} together with Remark \ref{Remark ii izole} are the main results of this section. \\

Recall that a subset $A$ of a ring $R$ is called a set of orthogonal idempotents if each member of $A$ is an idempotent and $e_{i}e_{k}=0$ for distinct elements $e_{i},e_{k}\in A$. In other words, $e_{i}e_{k}=\delta_{i,k}e_{i}$ for all $e_{i},e_{k}\in A$ where $\delta_{i,k}$ is the Kronecker delta. For example, if $e$ and $e'$ are idempotents of a ring $R$, then $\{e, 1-e\}$ and $\{e+e'-2ee', ee'\}$ are sets of orthogonal idempotents.

\begin{theorem}\label{Theorem  vii orthogonality} Lifting idempotents preserves the orthogonality in countable case.
\end{theorem}

{\bf Proof.} Let $\phi:R\rightarrow R'$ be a morphism of rings which lifts idempotents. We have to show that every countable set $\{e'_{1},e'_{2},e'_{3},\ldots\}$ of orthogonal idempotents of $R'$ can be lifted to a set $\{e_{1},e_{2},e_{3},\ldots\}$ of orthogonal idempotents of $R$. By the induction, it suffices to prove the assertion for each natural number $n$. So assume $\{e_{1},\ldots,e_{n}\}$ is a set of orthogonal idempotents of $R$ such that $\phi(e_{i})=e'_{i}$ for all $i=1,\ldots,n$. There exists an idempotent $e\in R$ such that $\phi(e)=e'_{n+1}$. Then clearly $e_{n+1}:=e(1-\sum\limits_{k=1}^{n}e_{k})$ is an idempotent,  $\phi(e_{n+1})=e'_{n+1}$ and $e_{n+1}e_{i}=0$ for all $i=1,\ldots,n$. $\Box$ \\

Let $X$ be a topological space and $S$ a topological ring. Then $\mathrm{C}(X,S)$, the set of all continuous functions $X\rightarrow S$, under the pointwise addition and multiplication is a ring. Indeed, if $f,g\in\mathrm{C}(X,S)$ then the induced map $h:X\mapsto S\times S$ given by $x\rightsquigarrow\big(f(x),g(x)\big)$ is continuous, because $h^{-1}(U\times V)=f^{-1}(U)\cap g^{-1}(V)$. The composition of $h$ with each of the additive and multiplicative operations of $S$ gives us $f+g$ and $f\cdot g$, respectively. Hence, $f+g$ and $f\cdot g$ are continuous.
The remaining ring conditions for $\mathrm{C}(X,S)$ are easily verified.
In particular, for a given ring $R$, setting $H_{0}(R):=\mathrm{C}\big(\Spec(R),\mathbb{Z}\big)$ where $\mathbb{Z}$ is equipped with the discrete topology. If  $\phi:R\rightarrow R'$ is a morphism of rings, then the map $H_{0}(\phi):H_{0}(R)\rightarrow H_{0}(R')$ given by $\psi\rightsquigarrow\psi\circ\phi^{\ast}$ is a morphism of rings.
In fact, $H_{0}(-)$ is a covariant functor from the category of commutative rings to itself. In the following result, we improve the formulation of \cite[Lemma 3.1(a)]{Bass} and provide a new proof for this lemma.

\begin{theorem}\label{Theorem lift h0} A morphism of rings $\phi:R\rightarrow R'$ lifts idempotents if and only if $H_{0}(\phi):H_{0}(R)\rightarrow H_{0}(R')$ is surjective. The map $H_{0}(\phi)$ is injective if and only if $\Ker\phi$ has no nonzero idempotents.
\end{theorem}

{\bf Proof.} Assume $\phi$ lifts idempotents. If $\psi\in H_{0}(R')$ then $\psi:\Spec(R')\rightarrow\mathbb{Z}$ is a continuous map. Thus for each $k\in\mathbb{Z}$, $\psi^{-1}(\{k\})$ is a clopen of $\Spec(R')$. Then by Theorem \ref{Theorem II}, there exists an idempotent $e'_{k}\in R'$ such that $\psi^{-1}(\{k\})=D(e'_{k})$. Using the quasi-compactness of the prime spectrum, then there exists a natural number $n\geqslant0$ such that $\Spec(R')=\bigcup\limits_{k=0}^{n}D(e'_{k})$. By hypothesis and Theorem \ref{Theorem  vii orthogonality}, there are orthogonal idempotents $e_{0},\ldots,e_{n}\in R$ such that $\phi(e_{k})=e'_{k}$ for all $k$. There exists also an idempotent $e_{n+1}\in R$ such that $D(e_{n+1})=\Spec(R)\setminus\bigcup\limits_{k=0}^{n}D(e_{k})$.
Now the map $\theta:\Spec(R)\rightarrow\mathbb{Z}$ defined by $\theta(\mathfrak{p})=k$ whenever $\mathfrak{p}\in D(e_{k})$ is well-defined, because the $D(e_{k})$'s are pairwise disjoint. It is also continuous and $\psi=\theta\circ\phi^{\ast}$. Thus $H_{0}(\phi)$ is surjective. Conversely, if $e'\in R'$ is an idempotent then the map $\psi:\Spec(R')\rightarrow\mathbb{Z}$ defined by $\psi(\mathfrak{p})$ is either $0$ or $1$, according to whether $e'\in\mathfrak{p}$ or $e'\notin\mathfrak{p}$, is continuous. Thus by hypothesis, there exists a continuous map $\theta:\Spec(R)\rightarrow\mathbb{Z}$ such that $\psi=\theta\circ\phi^{\ast}$. Then by Theorem \ref{Theorem II}, there exists an idempotent $e\in R$ such that $\theta^{-1}(\{1\})=D(e)$. It follows that $D\big(\phi(e)\big)=D(e')$. Thus by Theorem \ref{Theorem II}, $\phi(e)=e'$. Hence, $\phi$ lifts idempotents. Now we prove the second part of the assertion. Assume $H_{0}(\phi)$ is injective. Let $e\in \Ker\phi$ be an idempotent. Then the function $\psi:\Spec(R)\rightarrow\mathbb{Z}$ where $\psi(\mathfrak{p})$ is either $0$ or $1$, according as $e\in\mathfrak{p}$ or $e\notin\mathfrak{p}$, is continuous. Note that $\phi(e)=0$ and so $\psi\circ\phi^{\ast}=0$. Thus by hypothesis, $\psi=0$. It follows that $D(e)=\emptyset=D(0)$ and so $e=0$. Conversely, suppose there is a continuous function $\psi:\Spec(R)\rightarrow\mathbb{Z}$ such that $\psi\circ\phi^{\ast}=0$. We show that $\psi=0$. If not, then its image contains a nonzero integer $d$. By Theorem \ref{Theorem II}, there exists an idempotent $e\in R$ such that $\psi^{-1}(\{d\})=D(e)$. It follows that $D\big(\phi(e)\big)=\emptyset$ and so $\phi(e)=0$. Thus by hypothesis, $e=0$. But this is a contradiction, because $\psi^{-1}(\{d\})$ is nonempty. This completes the proof. $\Box$ \\

If $\phi:R\rightarrow R'$ is a morphism of rings, then the map $\mathcal{B}(\phi):\mathcal{B}(R)\rightarrow\mathcal{B}(R')$ is injective if and only if $\Ker\phi$ has no nonzero idempotents. In fact, we have the following result.

\begin{corollary} A morphism of rings $\phi:R\rightarrow R'$ strongly lifts idempotents if and only if $H_{0}(\phi):H_{0}(R)\rightarrow H_{0}(R')$ is bijective.
\end{corollary}

{\bf Proof.} If $\phi$ strongly lifts idempotents then by Theorem \ref{Theorem lift h0}, $H_{0}(\phi)$ is bijective.
Conversely, by Theorem \ref{Theorem lift h0}, $\phi$ lifts idempotents. Let $e,e'\in R$ be idempotents such that $\phi(e)=\phi(e')$. Then the idempotent $e+e'-2ee'\in\Ker\phi$. Thus $e+e'=2ee'$, because $\Ker\phi$ contains no nonzero idempotents. It follows that $e=ee'=e'$. $\Box$ \\

Remember that in a ring $R$, a nonzero idempotent $e\in R$ is called primitive if whenever $e'=ee'$ for some nonzero idempotent $e'\in R$, then $e'=e$. By Theorem \ref{Theorem II}, an idempotent $e\in R$ is primitive if and only if $\Clop\big(D(e)\big)\simeq\mathbb{Z}_{2}$, or equivalently, $D(e)$ is connected. In particular, in a nonzero ring $R$, then $1\in R$ is primitive if and only if $R$ has no nontrivial idempotents.

\begin{example}\label{Example I} We illustrate primitive idempotents with some basic examples:
\begin{enumerate}
\item Let $(R_{i})_{i\in I}$ be a family of rings and $R=\prod\limits_{i\in I}R_{i}$. Then the primitive idempotents of $R$ are precisely of the form $(e_{i})$ where $e_{k}\in R_{k}$ is a primitive idempotent for some $k\in I$ and $e_{i}=0$ for all $i\neq k$.
In particular, if each $R_{i}$ has no nontrivial idempotents, then the primitive idempotents of $R$ are precisely of the form $e_{k}=(\delta_{i,k})_{i\in I}$ for some $k\in I$.
\item Let $X$ be a set. As a particular case of (1), the primitive idempotents of the power set ring $\mathcal{P}(X)\simeq\prod\limits_{x\in X}\mathbb{Z}_{2}$ are precisely the singleton subsets of $X$.
\item In the ring $\mathbb{Z}_{30}=\mathbb{Z}/30\mathbb{Z}$, the idempotent $16$ is not primitive, since $6=6\times 16$ but $6\neq 16$.
\item By taking into account the Chinese remainder theorem and (1), then the idempotents and primitive idempotents of  $\mathbb{Z}_{m}=\mathbb{Z}/m\mathbb{Z}$ are implicitly computed. Here we give a formula to compute them in the explicit way. This formula can be applied for any Artinian ring (see Remark \ref{Example II}). In fact in Theorem \ref{Theorem IV} we provide a new and short proof to the Chinese remainder theorem. This formula is deduced from that proof. Consider the prime factorization $m=p^{c_{1}}_{1}\ldots p^{c_{n}}_{n}$ where the $p_{i}$'s are distinct prime numbers with $c_{i}\geqslant1$. Then the ring $\mathbb{Z}_{m}$ has $2^{n}$ idempotents and $($modulo $m$$)$ they are precisely of the form $\sum\limits_{k=1}^{n}h_{k}\epsilon_{k}$ where $\epsilon_{k}\in\{0,1\}$ and $h_{k}\in(\prod\limits_{\substack{i=1,\\
i\neq k}}^{n}p^{c_{i}}_{i})\mathbb{Z}$ such that $h_{k}-1\in p^{c_{k}}_{k}\mathbb{Z}$. Also $\{h_{1},\ldots,h_{n}\}$ is the set of primitive idempotents of $\mathbb{Z}_{m}$.
\item As a specific example of (5), we have $30=2\times 3\times 5$ so the idempotents of $\mathbb{Z}_{30}$ are precisely of the form $15\epsilon_{1}+10\epsilon_{2}+6\epsilon_{3}$. Thus $\{0,1,6,10,15,16,21,25\}$ is the set of idempotents and $\{6,10,15\}$ is the set of primitive idempotents of $\mathbb{Z}_{30}$.
\item As another specific example of (5), we have $765=3^{2}\times 5\times17$. Hence,
the idempotents of $\mathbb{Z}_{765}$ are precisely of the form
$306\epsilon_{1}-135\epsilon_{2}-170\epsilon_{3}$. Therefore $\{0,1,136,171,306,460,595,630\}$ is the set of idempotents and $\{306, 595, 630\}$ is the set of primitive
idempotents of $\mathbb{Z}_{765}$.
\item The primitive idempotents need not pass to extension of rings. For example, $1\in\mathbb{Z}_{2}$ is primitive, but it is not primitive in the strict Boolean extensions of $\mathbb{Z}_{2}$.
\end{enumerate}
\end{example}

\begin{theorem}\label{Theorem x isolated points} The primitive idempotents of a zero dimensional ring $R$ are in bijection with the isolated points of $\Spec(R)$.
\end{theorem}

{\bf Proof.} Clearly the primitive idempotents of a ring $R$ and its Boolean ring $\mathcal{B}(R)$ are the same. Also by Lemma \ref{Lemma vi zero dim}(iv), the isolated points of $\Spec(R)$ and $\Spec\big(\mathcal{B}(R)\big)$ are the same.
Therefore without loss of generality, we may assume that $R$ is a Boolean ring. Then it can be shown that an element $e\in R$ is primitive if and only if $Re=\{0,e\}$ is nonzero, or equivalently, $\mathfrak{p}_{e}=R(1-e)$ is a prime ideal of $R$. Now we show that the map $e\rightsquigarrow \{\mathfrak{p}_{e}\}$ is a one-to-one correspondence from the set of primitive idempotents of $R$ onto the set of isolated points of $\Spec(R)$. If $e$ is a primitive of $R$ then $\{\mathfrak{p}_{e}\}=D(e)$, because if $\mathfrak{q}\in D(e)$ then $1-e\in\mathfrak{q}$ thus $\mathfrak{p}_{e}=R(1-e)\subseteq\mathfrak{q}$ and so $\mathfrak{p}_{e}=\mathfrak{q}$. Hence the above map is well defined. It is also injective. Finally, if $\mathfrak{q}$ is an isolated point of $\Spec(R)$ then $\{\mathfrak{q}\}=D(e)$ for some $e\in R$. We show that $e$ is primitive. Clearly it is nonzero. If $ee'=e'$ for some nonzero $e'\in R$, then $D(e)=D(e')$ and so by Theorem \ref{Theorem II}, $e=e'$. $\Box$ \\

The above result does not hold for rings with positive Krull dimension. For example, $1\in\mathbb{Z}$ is primitive but $\Spec(\mathbb{Z})$ has no isolated points. Also note that the converse of Theorem \ref{Theorem x isolated points} does not hold. In other words, there are rings $R$ with $\dim(R)>0$ such that the set of primitive idempotents of $R$ corresponds bijectively with the set of isolated points of $\Spec(R)$. For example, take $R=\mathcal{P}(X)\times\mathbb{Z}$ where $X$ is an infinite set.

\begin{remark}\label{Remark ii izole} If $R$ is a finite Boolean ring, then $R\simeq\mathcal{P}(X)$ where $X=\Spec(R)$. So $|R|=2^{n}$ where $n=|\Spec(R)|$ is the number of primitive elements of $R$. But it is important to notice that if the set of primitive idempotents of a Boolean ring $R$ is finite, then $R$ is not necessarily finite. For example, let $2^{\mathbb{N}}\simeq\mathcal{P}(\mathbb{N})$ be the countably infinite products of the discrete space $\{0,1\}$ where $\mathbb{N}=\{0,1,2,\ldots\}$ is the set of natural numbers. The set $2^{\mathbb{N}}$ carries the product topology and so by Tychonoff theorem, it is compact (quasi-compact and Hausdorff). It is also totally disconnected, since every product of totally disconnected spaces is totally disconnected. So by Stone duality (see e.g. \cite[Corollary 5.5]{Tarizadeh 2}), there exists a Boolean ring $R$ such that $\Spec(R)$ is homeomorphic to $2^{\mathbb{N}}$. The space $2^{\mathbb{N}}$ has no isolated points. Because if $s=(s_{i})_{i\geqslant0}$ is an isolated point of $2^{\mathbb{N}}$ then there is a base open $U=\prod\limits_{i\geqslant0}U_{i}$ of $2^{\mathbb{N}}$ such that $\{s\}=U$ where each $U_{i}=\{0,1\}$ except for finitely many $U_{i_{1}}=\{s_{i_{1}}\},\ldots,U_{i_{n}}=\{s_{i_{n}}\}$. Now consider the sequence $t=(t_{i})_{i\geqslant0}\in2^{\mathbb{N}}$ where $t_{i}=1-s_{i}$ for all $i\in\mathbb{N}\setminus\{i_{1},\ldots,i_{n}\}$ and $t_{i_{k}}=s_{i_{k}}$ for all $k\in\{1,\ldots,n\}$. Then clearly $t\neq s$ and $t\in U$. This is a contradiction. Therefore by Theorem \ref{Theorem x isolated points}, $R$ has no primitive idempotents. But $R$ is infinite, because $\Spec(R)\simeq2^{\mathbb{N}}$ is infinite. Using this Boolean ring $R$, then for each $n\geqslant1$, the infinite Boolean ring $R\times\prod\limits_{i=1}^{n}\mathbb{Z}_{2}$ has precisely $n$ primitive idempotents. \\
\end{remark}

The above observation leads us to the following result.

\begin{theorem}\label{theorem xi primitive} If a ring $R$ as an ideal is generated by a set of its primitive idempotents, then $R$ has finitely many idempotents.
\end{theorem}

{\bf Proof.} If $e$ and $e'$ are distinct primitive idempotents of a ring $R$, then $ee'=0$. Assume $ee'\neq0$, we have $(ee')e'=ee'$ hence $ee'=e'$, since $e'$ is primitive. But $e$ is also primitive, so $e=e'$, a contradiction. Therefore, every set of (distinct) primitive idempotents is orthogonal. By hypothesis, there exists a finite set $\{e_{1},\ldots,e_{n}\}$ of (distinct) primitive idempotents of $R$ such that $\sum\limits_{i=1}^{n}r_{i}e_{i}=1$. Setting $e=\sum\limits_{k=1}^{n}e_{k}$. Then
$e_{1}\oplus\cdots\oplus e_{n}=e=\sum\limits_{i=1}^{n}r_{i}ee_{i}=
\sum\limits_{i=1}^{n}r_{i}e_{i}=1$. If $\mathfrak{p}\in X=\Spec\big(\mathcal{B}(R)\big)$ then $e_{i}\notin\mathfrak{p}$ for some $i$. It follows that $1-e_{i}\in\mathfrak{p}$ and so $(1-e_{i})\subseteq\mathfrak{p}$. But $(1-e_{i})$ is a prime ideal of $\mathcal{B}(R)$ and so $(1-e_{i})=\mathfrak{p}$. Therefore $X$ is finite. Hence $\mathcal{B}(R)$ is finite, since by Theorem \ref{Theorem II}, $\mathcal{B}(R)\simeq\Clop(X)$. $\Box$

\begin{proposition} If each regular ideal of a ring $R$ is finitely generated, then there exists a finite set $\{e_{1},\ldots,e_{n}\}$ of orthogonal idempotents of $R$ such that $\sum\limits_{i=1}^{n}e_{i}=1$ and $|\mathcal{B}(R)|=2^{n}$.
\end{proposition}

{\bf Proof.} If $C$ is a connected component of $\Spec(R)$, then by  hypothesis, there exists an idempotent $e\in R$ such that $C=D(e)$. Now using the quasi-compactness of $\Spec(R)$, then there exists a finite set $\{e_{1},\ldots,e_{n}\}$ of orthogonal idempotents of $R$ such that $\Spec(R)=\bigcup\limits_{i=1}^{n}D(e_{i})$ and $n$ is the number of connected components of $\Spec(R)$. Thus by Lemma \ref{Theorem ix finiteness}, $|\mathcal{B}(R)|=2^{n}$. We also have $R=(e_{1},\ldots,e_{n})$.
Thus we may write $1=\sum\limits_{k=1}^{n}r_{k}e_{k}$. Setting $e=\sum\limits_{i=1}^{n}e_{i}$. Then $e=\sum\limits_{k=1}^{n}r_{k}ee_{k}=1$. $\Box$ \\

Every Noetherian ring or more generally every semi-Noetherian ring in the sense of \cite[Definition 6.1]{Tarizadeh 5} satisfies in the hypothesis of the above proposition.

\begin{theorem}\label{Theorem xii TD} Let $\{\mathfrak{p}_{1},\ldots,\mathfrak{p}_{n}\}$ be a finite set of distinct prime ideals of a zero dimensional ring $R$ where $n\geqslant1$. Then there exists a set $\{e_{1},\ldots,e_{n}\}$ of orthogonal idempotents of $R$ such that each $\mathfrak{p}_{i}\in D(e_{i})$ and $\sum\limits_{i=1}^{n}e_{i}=1$.
\end{theorem}

{\bf Proof.} If $n=1$ then it suffices to take $e_{1}=1$. Thus we may assume $n\geqslant2$. The space $X=\Spec(R)$ is compact and totally disconnected. So the set $\Clop(X)$ forms a base for the opens of $X$. Therefore if $i\neq k$ then we may choose disjoint clopen neighborhoods $\mathfrak{p}_{i}\in U_{i,k}$ and $\mathfrak{p}_{k}\in U_{k,i}$ of $X$. Setting $V_{i}=\bigcap\limits_{\substack{k=1,\\
k\neq i}}^{n}U_{i,k}$ for $i=1,\ldots,n-1$ and $V_{n}=X\setminus\bigcup\limits_{i=1}^{n-1}V_{i}$. Then clearly $V_{1},\ldots,V_{n}$ are pairwise disjoint clopens of $X$ such that $\mathfrak{p}_{i}\in V_{i}$ for all $i$. By Theorem \ref{Theorem II}, there exists a set $\{e_{1},\ldots,e_{n}\}$ of orthogonal idempotents of $R$ such that each $V_{i}=D(e_{i})$. Clearly the ideal $(e_{1},\ldots,e_{n})$ is the whole ring $R$. It follows that $\sum\limits_{i=1}^{n}e_{i}=1$. $\Box$

\begin{remark}\label{Remark iii reverse} The converse of Theorem \ref{Theorem xii TD} also holds. In fact, any finite set of prime ideals of a ring $R$ makes a partition of its prime spectrum into clopen parts with distinct elements of this set inside distinct parts if and only if $R$ is zero dimensional. Indeed, the implication $``\Leftarrow"$ is just Theorem \ref{Theorem xii TD}. Conversely, if $\mathfrak{p}\subset\mathfrak{q}$ are distinct primes of $R$, then by hypothesis and Theorem \ref{Theorem II}, there is an idempotent $e\in\mathfrak{q}\setminus\mathfrak{p}$. It follows that $1-e\in\mathfrak{p}$ which is a contradiction. Hence, $\dim(R)=0$.
\end{remark}

\begin{example} If $\{\mathfrak{p}_{1},\ldots,\mathfrak{p}_{n}\}$ is a finite set of distinct minimal primes of a ring $R$, then by the prime avoidance lemma (see e.g. \cite[Theorem 2.2]{Tarizadeh-Chen}), $\Spec(S^{-1}R)=\{S^{-1}\mathfrak{p}_{1},\ldots,S^{-1}\mathfrak{p}_{n}\}$ where $S=R\setminus\bigcup\limits_{i=1}^{n}\mathfrak{p}_{i}$. Hence, $S^{-1}R$ is zero dimensional. If the canonical ring map $R\rightarrow S^{-1}R$ lifts idempotents, then using Theorems \ref{Theorem  vii orthogonality} and \ref{Theorem xii TD}, we may find orthogonal idempotents $e_{1},\ldots,e_{n}\in R$ and some $s\in S$ such that each $\mathfrak{p}_{i}\in D(e_{i})$ and $s(1-\sum\limits_{i=1}^{n}e_{i})=0$.
Dually, let $\{\mathfrak{m}_{1},\ldots,\mathfrak{m}_{d}\}$ be a finite set of distinct maximal ideals of a ring $R$. If $I=\bigcap\limits_{i=1}^{d}\mathfrak{m}_{i}$ is a lifting ideal, then there exist orthogonal idempotents $e'_{1},\ldots,e'_{d}\in R$ such that each $\mathfrak{m}_{k}\in D(e'_{k})$ and $1-\sum\limits_{k=1}^{d}e'_{k}\in I$. \end{example}

\begin{proposition} Let $e$ be a nonzero idempotent of a ring $R$. Then $e$ is primitive if and only if it has no decomposition as $e=e_{1}+e_{2}$ where $e_{1},e_{2}\in R$ are nonzero orthogonal idempotents.
\end{proposition}

{\bf Proof.} Assume $e$ is primitive. If we have a decomposition as the above, then $(1-e)e_{1}=0$. If $e_{1}\neq0$ then $e_{1}=e$ and so $e_{2}=0$. Conversely, if $(1-e)e'=0$ for some nonzero idempotent $e'\in R$. Then we may write $e=e'+(1-e')e$. Thus $(1-e')e=0$ and so $e=e'$. Hence, $e$ is primitive. $\Box$

\begin{remark} Let $\phi:R\rightarrow R'$ be a morphism of rings whose kernel is contained in the Jacobson radical of $R$. If $e\in R$ is an idempotent such that $\phi(e)$ is primitive, then $e$ is also primitive. Because if $e'=ee'$ for some nonzero idempotent $e'\in R$. Then $\phi(e')\neq0$ and so $e-e'\in\mathfrak{J}$. Thus by Lemma \ref{Lemma II}, $e=e'$. The above hypothesis on the ring map is also crucial to characterize the projectivity of finitely generated flat modules, see \cite[Theorem 3.3]{Abolfazl extracta}.
\end{remark}

\section{Appendix}

Here we give a new and significantly short proof to the Chinese remainder theorem. This proof is also constructive and provides a formula to compute the idempotents of an Artinian ring (especially $\mathbb{Z}/m\mathbb{Z}$) in the explicit way, see Example \ref{Example I}(4) and Remark \ref{Example II}.

\begin{theorem}$($Chinese Remainder Theorem$)$\label{Theorem IV} Let $I_{1},\ldots,I_{n}$ be finitely many proper ideals of a ring $R$ which are pairwise coprime $($i.e., $I_{i}+I_{k}=R)$. Then the canonical ring map $\pi:R\rightarrow R/I_{1}\times\cdots\times R/I_{n}$ given by $f\rightsquigarrow(f+I_{1},\ldots,f+I_{n})$ is surjective and its kernel is $\prod\limits_{i=1}^{n}I_{i}$.
\end{theorem}

\begin{proof} Clearly for each $k$, we have $I_{k}+\bigcap\limits_{\substack{i=1,\\
i\neq k}}^{n}I_{i}=R$. Thus there exist $g_{k}\in I_{k}$ and $h_{k}\in\bigcap\limits_{\substack{i=1,\\
i\neq k}}^{n}I_{i}$ such that $1=g_{k}+h_{k}$. If $(f_{1},\ldots,f_{n})$ is an $n$-tuple of elements of $R$ then for each $k$, one has $f-f_{k}=f-f_{k}g_{k}-f_{k}h_{k}=\sum\limits_{\substack{i=1,\\
i\neq k}}^{n}f_{i}h_{i}-f_{k}g_{k}\in I_{k}$ where $f=\sum\limits_{i=1}^{n}f_{i}h_{i}$. Therefore $\pi(f)=(f_{1}+I_{1},\ldots,f_{n}+I_{n})$. Also clearly $\Ker\pi=\bigcap\limits_{i=1}^{n}I_{i}=\prod\limits_{i=1}^{n}I_{i}$.
\end{proof}

\begin{remark}\label{Example II} Let $R$ be an Artinian ring. Then every prime ideal of $R$ is a maximal ideal and it has finitely many maximal ideals. Let $\{\mathfrak{m}_{1},\ldots,\mathfrak{m}_{n}\}$ be the distinct maximal ideals of $R$. The Jacobson radical of $R$ is nilpotent. Let $N$ be the least natural number such that $\mathfrak{J}^{N}=0$. The map $R\rightarrow R/\mathfrak{m}^{N}_{1}
\times\cdots\times R/\mathfrak{m}^{N}_{n}$ given by $f\rightsquigarrow
(f+\mathfrak{m}^{N}_{1},\ldots,f+\mathfrak{m}^{N}_{n})$ is an isomorphism of rings. Using this, then the proof of Theorem \ref{Theorem IV} shows that $f\in R$ is an idempotent if and only if $f=\sum\limits_{i=1}^{n}h_{k}\epsilon_{k}$ where $\epsilon_{k}\in\{0,1\}$ and $h_{k}\in\prod\limits_{\substack{i=1,\\
i\neq k}}^{n}\mathfrak{m}^{N}_{i}$ such that $h_{k}-1\in\mathfrak{m}^{N}_{k}$. Moreover $\{h_{1},\ldots,h_{n}\}$ is the set of primitive idempotents of $R$. In particular, $R$ has $2^{n}$ idempotents and $n$ primitive idempotents. \\
\end{remark}

\textbf{Acknowledgements.} We express our deep thanks to Professor Pierre Deligne for the very fruitful discussions that we had with him during the writing of the present article. We would also like to give sincere thanks to the referee for very careful reading of the article, and for his/her very valuable suggestions and comments.

\end{document}